\documentclass[10pt]{amsart}
\usepackage{amssymb}
\usepackage[curve,matrix,arrow]{xy}
\numberwithin{equation}{section}
\newtheorem{thm}{Theorem}[section]
\newtheorem{prop}[thm]{Proposition}
\newtheorem{fact}[thm]{Fact}
\newtheorem{cor}[thm]{Corollary}
\newcommand\Hom{\operatorname{Hom}}

\newcommand\im{\operatorname{im}}
\newcommand\comp{\smash{\lower-.1ex\hbox{\scriptsize$\, \circ\, $}}}
\newcommand{\Z}{{\mathbb Z}}
\newcommand{\Q}{{\mathbb Q}}
\title{Splitting off Rational Parts in Homotopy Types}
\author[Iwase]{%
Norio Iwase${}^\dagger$}
\email{iwase@math.kyushu-u.ac.jp}
\author[Oda]{%
Nobuyuki Oda${}^{\dagger\dagger}$}
\email{odanobu@cis.fukuoka-u.ac.jp}
\address[Iwase]{%
Faculty of Mathematics,
Kyushu University,
Fukuoka 810-8560, Japan}
\address[Oda]{%
Department of Applied Mathematics, 
Fukuoka University, 
Fukuoka 814-0180, Japan}
\date\today
\thanks{${}^\dagger$supported by the Grant-in-Aids for Scientific Research \#14654016 from the Ministry of Education, Culture, Sports, Science and Technology, Japan.}
\thanks{${}^{\dagger\dagger}$supported by the Grant-in-Aids for Scientific Research \#15340025 from the Japan Society for the Promotion of Science.}
\begin{document} 
\begin{abstract}
It is known algebraically that any abelian group is a direct sum of a divisible group and a reduced group (See Theorem 21.3 of \cite{Fuchs:abelian-group}).
In this paper, conditions to split off rational parts in homotopy types from a given space are studied in terms of a variant of Hurewicz map, say $\overline{\rho} : [S_{\Q}^{n},X] \to H_n(X;\Z)$ and generalized Gottlieb groups.
This yields decomposition theorems on rational homotopy types of Hopf spaces, $T$-spaces and Gottlieb spaces, which has been known in various situations, especially for spaces with finiteness conditions.
\end{abstract}
\maketitle

\section*{Introduction} 

The Gottlieb group is introduced by Gottlieb \cite{Gottlieb65a,Gottlieb69a} and the generalized Gottlieb set is introduced by Varadarajan \cite{Varadarajan69}. 
Dula and  Gottlieb obtained a general result on splitting a Hopf space off from a fibration as Theorem 1.3 of \cite{DuGo90}.

In this paper, we work in the category of spaces having homotopy types of CW complexes with base points and pointed continuous maps.
A relation $f \sim g$ indicates a pointed homotopy relation of maps $f$ and $g$ and a relation $X \simeq Y$ indicates a homotopy equivalence relation of spaces $X$ and $Y$.
We also denote by $[X,Y]$ the set of pointed homotopy classes of maps from $X$ to $Y$. 

We adopt some more conventional notations:
$X_{\Q}$ stands for the rationalization of a space $X$, 
$K(\pi,n)$ for the Eilenberg-Mac\,Lane space of type $(\pi,n)$, $G(V,X)$ for the generalized Gottlieb subset of $[V,X]$ and $H_n(X)$ for $H_n(X;\Z)$.
We introduce a variant of Hurewicz map 
$\overline{\rho} : [S_{\Q}^{n},X] \to H_n(X)$ by $\overline{\rho} (\alpha) = \alpha_*([S^n]{\otimes}1)$ for $\alpha \in [S_{\Q}^{n},X] $, where $\alpha_*$ is the homomorphism given by $\alpha_{*} : H_n(S^{n}){\otimes}{\Q} = H_n(S_{\Q}^{n}) \to H_n(X)$.
Our main result is described as follows:
 
\par\medskip\noindent 
{\bf Theorem \ref{th:rathopfdec}.  \it 
Let $R=\bigoplus_{\lambda \in \Lambda}\Q\,$ be a $\Q$-vector space of dimension $\#\Lambda \leq \infty$.
Let $X$ be $0$-connected and $R \subset \overline{\rho}(G(S_{\Q}^{n},X)) \subseteq H_n(X)$, $n\geq 2$.
Then $X$ decomposes as 
$$
X 
\simeq 
Y \times K(R,n).
$$ 
}
\par
Theorem \ref{th:rathopfdec} gives unified proof to the splitting phenomena on rational G-space, $T$-space and Hopf space without assuming any finiteness conditions, which are proved under various situations by a number of authors: 
Scheerer \cite{Scheerer85} obtained decomposition theorems 
of rational Hopf spaces without assuming the finite type assumptions.  
Oprea \cite{Oprea86} obtained decomposition theorems by using minimal model method in rational homotopy theory.
Aguad\'{e} \cite{Aguade87} obtained a decomposition theorems on rational $T$-spaces of finite type.

\section{Preliminaries}\label{Secrational}

We regard the one point union $X {\vee} Y$ of spaces $X$ and $Y$ as a subspace $X{\times}{\ast} \cup {\ast}{\times}Y$ of the product space $X {\times} Y$ with the inclusion map $j : X {\vee} Y \rightarrow X {\times} Y$. 
For any collection of a finitely or infinitely many spaces $X_{\lambda}$ ($\lambda\in\Lambda$), we denote the {\it wedge sum} (or one point union) by  $\bigvee_{\lambda\in\Lambda} X_{\lambda}$ and the {\it direct sum} (or weak product) by 
$\bigoplus_{\lambda\in\Lambda}X_{\lambda} = 
\left\{\left.
(x_{\lambda}) \in \prod_{\lambda\in\Lambda}X_{\lambda}
 \  \right| \ 
x_{\lambda}=\ast \ \ \text{except for finitely many $\lambda$}
\right\}$.
Then we have $\bigvee_{\lambda\in\Lambda} X_{\lambda} \subset \bigoplus_{\lambda\in\Lambda}X_{\lambda}$, where $\bigoplus_{\lambda\in\Lambda}X_{\lambda}$ is a dense subset of the product space $\prod_{\lambda\in\Lambda}X_{\lambda}$ and has the weak topology with respect to finite products of $X_{\lambda}$'s. 

We can apply rationalization or $\Q$-localization to any {\it $0$-connected nilpotent spaces} (cf. \cite{BoKa72}, \cite{HiMiRo75} or \cite{MiNiTo71}). 
The  rationalization $\ell_{\Q} : X \to X_\Q$, or simply $X_{\Q}$ does exist for such spaces $X$ such that $\ell_{\Q}$ induces the following isomorphisms:
$$
\pi_n(X_\Q) \cong \pi_n(X) \otimes \Q  \mbox{\ \ \ and \ \ \ } 
H_n(X_\Q) \cong H_n(X) \otimes \Q$$ 
for any integer $n\geq 1$.
Moreover the universality of rationalization yields a bijection
$$
\ell^{\ast}_{\Q} : [X_\Q , Y_\Q] \cong [X , Y_\Q]
$$
for any such spaces $X$ and $Y$.  
The rationalization enjoys the following fact.

\begin{fact}\label{fact:rational} 
\begin{enumerate}
\item\label{fact.1}
$S_{\Q}^{2m+1} \simeq K(\Q, 2m+1)$  for any integer $m\geq 0$.
\item\label{fact.2}
$\Omega(S_{\Q}^{2m+1}) \simeq ({\Omega}S^{2m+1})_{\Q} \simeq K(\Q, 2m)$  for any integer $m\geq 1$. 
\item\label{fact.3}
$(X_\infty)_\Q \simeq (X_\Q)_\infty$ for $X$ a  $0$-connected  
nilpotent space of finite type. 
\end{enumerate}
\end{fact}
\begin{proof}
(\ref{fact.1}) and (\ref{fact.2}) are well-known.
We give here a brief explanation for (\ref{fact.3}):
The suspension functor $\Sigma$ and the loop functor $\Omega$ enjoys the properties 
\ $\Sigma(X_\Q) \simeq (\Sigma X)_\Q$  \ for any $0$-connected space $X$ and  \  $\Omega(X_\Q) \simeq (\Omega X)_\Q$  \  
for any $1$-connected space $X$. 
Let $X_\infty$ be the James reduced product space of a $0$-connected space $X$ of finite type, so that $X_\infty \simeq \Omega (\Sigma X)$ by James \cite{James55}.
Then it follows that 
$(X_\infty)_\Q \simeq (\Omega (\Sigma X))_\Q \simeq  \Omega (\Sigma (X_\Q)) \simeq (X_\Q)_{\infty}$.
\end{proof}

We state two propositions to be used in the proof of the main theorem.
\begin{prop} \label{pr:jamesredprodspext} 
Let $X$ be a $0$-connected space of finite type and $f : X \to Y$ a map.
If $f \in G(X,Y)$, then there is an extension $\overline{f} : X_\infty \to Y$ of $f$ such that $\overline{f} \in G(X_\infty , Y)$.  
\end{prop} 

\begin{proof} 
We may assume that there is a map $\mu : Y \times X \to Y$ such that $\mu | Y \times \{ * \} = 1_Y:  Y \to Y$ and $\mu | \{ * \} \times X = f : X  \to Y$.
We put $\mu_1 = \mu$ and, for any $n$ we define  
$$
\mu_n = \mu \comp (\mu_{n-1} \times 1_X) : 
Y \times X^{n} = (Y \times X^{n-1}) \times X \to Y \times X 
\to Y
$$
by induction on $n$.
Then we observe that $\mu_n$ factors through $Y \times X^{n} \to Y \times X_n$.
\end{proof}

\begin{prop} \label{pr:infpairingext} 
Let $\alpha_\lambda : X_\lambda \to Z$ be a map for any $\lambda \in \Lambda$.  
If $\alpha_\lambda \in G(X_\lambda,Z)$ for each $\lambda \in \Lambda$, 
then the map $\alpha : \bigvee_{\lambda \in \Lambda} X_\lambda  \to Z$ defined by 
$\alpha |  X_\lambda = \alpha_\lambda  : X_\lambda \to Z$  can be extended to a map 
$\overline{\alpha} : \bigoplus_{\lambda \in \Lambda} X_\lambda  \to Z$ with 
$\overline{\alpha} \in G( \bigoplus_{\lambda \in \Lambda}X_\lambda,Z)$. 
\end{prop} 
\begin{proof} 
Since each $X_\lambda$ has a homotopy type of a CW complex, we may assume that there is a map $\mu_{\lambda} : Z \times X_{\lambda} \to Z$ such that $\mu_{\lambda} | \{ * \} \times X_{\lambda}  = \alpha_{\lambda} : X_{\lambda}  \to Z$ and $\mu_{\lambda} | Z \times \{ * \}   = 1_{Z} : Z  \to Z$ for each  $\lambda \in \Lambda$.  
For any $n$ and $\lambda_1, \lambda_2, \cdots, \lambda_n$, we define 
$$
\mu_{\lambda_{1},\cdots,\lambda_{n}} = \mu_{\lambda_n} \circ (\mu_{\lambda_{1},\cdots,\lambda_{n-1}} \times 1_{X_{\lambda_n}}) 
: Z \times ({X_{\lambda_1}} \times \cdots \times {X_{\lambda_{n-1}}} \times {X_{\lambda_n}}) \to Z
$$
by induction on $n$.
For any index set $\Lambda$, we assume that $\Lambda$ is totally-ordered.
Then we can easily observe that the collection of maps $\mu_{\lambda_{1},\cdots,\lambda_{n}}$ defines a pairing $\mu : Z \times (\bigoplus_{\lambda \in \Lambda}X_\lambda ) \to Z$ with axes $(1_Z , \overline{\alpha})$ (cf. \cite{Oda92}).
\end{proof}  

\section{Proof of the main result}\label{SecGottliebgroup}

\begin{prop}\label{pr:albeid} 
Let $P$ be an idempotent endomorphism of $H_n(X)$, $n \geq 2$. 
Suppose that $R = \im{P} \subseteq H_n(X)$ is a rational vector space and is in $\im{\overline{\rho}}$.
Then we have maps $\alpha : S^{n}(R) \to X$ and $\beta : X \to K(R,n)$ such that 
\begin{align*}&
\beta \comp \alpha \sim \iota^{n}_R : S^{n}(R) \to K(R,n), \ \text{and}
\\&
P = \alpha_{\ast} \comp (\iota^{n}_{R\,\ast})^{-1} \comp \beta_{\ast} : H_n(X) \to H_n(K(R,n)) \overset{\cong}{\leftarrow} H_n(S^{n}(R)) \to H_n(X),
\end{align*}
where $S^{n}(R)$ denotes the Moore space of type $(R,n)$ and $\iota^n_{R}$ corresponds to the identity element in $\Hom(R,R) = \Hom(\pi_n(S^{n}(R)),\pi_n(K(R,n))) \cong [S^{n}(R),K(R,n)]$. 
\end{prop}

\begin{proof}
Let $\{\overline{\rho}(\alpha_{\lambda}) \,\vert\, \lambda \in \Lambda\}$ be a basis of $R = \im{P}$, and hence $R \cong \bigoplus_{\lambda\in\Lambda} \Q$.
Since $S^{n}(R) = \bigvee_{\lambda\in\Lambda} S^{n}_{\Q}$,  we define $\alpha : S^{n}(R) \to X$ by its restrictions to all factors: 
$$
\alpha\vert_{S^{n}_{\Q}} = \alpha_{\lambda} : S^{n}_{\Q} \to X.
$$
Since $\alpha_{\ast}$ is an isomorphism onto $R \subseteq H_n(X)$, we have its inverse $\phi : R \to H_n(S^{n}_{\Q})$ so that $\phi \comp \alpha_* = {\rm id}_{H_n(S^{n}_{\Q})}$ and $\alpha_* \comp \phi = {\rm id}_{R}$.
Now we define a homomorphism $s : H_n(X) \to \im P \cong H_n(S^{n}_{\Q})$ by $s = \phi \comp P$: Since $\im{\alpha_{\ast}}$ is in the image of an idempotent endomorphism $P$, we have $s \comp \alpha_* = \phi \comp P \comp \alpha_* = \phi \comp \alpha_* = {\rm id}$.
Also we have $\alpha_* \comp s = \alpha_* \comp \phi \comp P = P$.
Thus $s$ satisfies the following formulae:
\begin{align*}&
s \comp \alpha_* = {\rm id} : H_n(S^{n}(R)) \to H_n(S^{n}(R)),
\\&
\alpha_* \comp s = P : H_n(X) \to H_n(X). 
\end{align*}

Let us recall that $\alpha$ induces the following commutative diagram:

\begin{equation}\label{cd:hurewicz}
\xymatrix{
[X,K(R,n)]
        \ar[d]_{\alpha^*}
        \ar[r]^{\hspace*{-10ex}\Psi'}_{\hspace*{-10ex}\cong}
&
{\rm Hom}(H_n(X),H_n(K(R,n)))
        \ar[d]^{(\alpha_*)^*}
\\
[S^{n}(R),K(R,n)]%
        \ar[r]^{\hspace*{-10ex}\Psi}_{\hspace*{-10ex}\cong}
&
{\rm Hom}(H_n(S^{n}(R)),H_n(K(R,n)))%
,
}
\end{equation}
where $\Psi$ and $\Psi'$ are homomorphisms defined by  taking the $n$-th homology groups, and are isomorphisms by the universal coefficient theorem.  
Since $\Psi'$ is an isomorphism, we define $\beta$ to be the unique element ${\Psi'}^{-1}(\iota^{n}_{R\ast} \comp s)$ so that $\beta_{\ast} = \iota^n_{R\ast} \comp s$.

Firstly by $P = \alpha_{\ast} \comp s$, we have $P = \alpha_{\ast} \comp s = \alpha_{\ast} \comp (\iota^{n}_{\Q\,\ast})^{-1} \comp \beta_{\ast}$.

Next we show $\beta \comp \alpha \sim \iota^{n}_{\Q}$.
By the commutativity of the diagram (\ref{cd:hurewicz}), we have 
\begin{align*}
\Psi(\alpha^*(\beta)) &= (\alpha_*)^{*}\comp\Psi'(\beta) 
= (\alpha_*)^*(\iota^{n}_{\Q\,\ast} \comp s)
=\iota^{n}_{\Q\,\ast} \comp  s \comp \alpha_* 
= \iota^{n}_{\Q\,\ast}
= \Psi(\iota^{n}_{\Q}).  
\end{align*}
Since $\Psi$ is an isomorphism, we also have 
$\beta \comp \alpha = \alpha^*(\beta) \sim \iota^{n}_{\Q}$. 
\end{proof}

Let us recall that $G(S_{\Q}^{n},X) \subset [S_{\Q}^{n},X] \xrightarrow{\overline{\rho}} H_n(X)$.  
In the following theorem, we do {\it not} assume that $X$ is rationalized nor that $X$ is $(n{-}1)$-connected.

\begin{thm}\label{th:rathopfdec}
Let $R=\bigoplus_{\lambda \in \Lambda}\Q\,$ be a $\Q$-vector space of dimension $\#\Lambda \leq \infty$.
Let $X$ be $0$-connected and $R \subset \overline{\rho}(G(S_{\Q}^{n},X)) \subseteq H_n(X)$, $n\geq 2$.
Then $X$ decomposes as 
$$
X 
\simeq 
Y \times K(R,n).
$$ 
\end{thm} 
\begin{proof} 
Since a divisible submodule $R$ is a direct summand of $H_n(X)$, there is an idempotent endomorphism $P : H_n(X) \to H_n(X)$ with $\im{P} = R$.
We fix a basis of $R$ as  $\{\overline{\rho}(\alpha_{\lambda}) \,\vert\, \alpha_{\lambda} \in G(S^{n}_{\Q},X), \lambda\in\Lambda\}$.

By Proposition \ref{pr:albeid}, there are maps $\alpha : S^{n}(R) \to X$, $\beta : X \to K(R,n)$ such that 
\begin{align*}&
\beta \comp \alpha \sim \iota^{n}_R : S^{n}(R) \to K(R,n), 
\\&
P = \alpha_{\ast} \comp (\iota^{n}_{R\,\ast})^{-1} \comp \beta_{\ast} : H_n(X) \to H_n(K(R,n)) \overset{\cong}{\leftarrow} H_n(S^{n}(R)) \to H_n(X).
\end{align*}
Then we extend the map $\alpha$ onto $K(R,n) \supseteq S^n(R)$ as $\overline{\alpha} : K(R,n) \to X$ by dividing our arguments in two cases: 
\par\noindent(Case 1)
$n$ is an odd positive integer $> 1$, namely,  $n = 2m+1$ for some $m \geq 1$.
Then we have $K(\Q,2m+1) \simeq S^{2m+1}_{\Q}$, and hence we have nothing to do.
\par\noindent(Case 2)
$n$ is an even positive integer, namely, $n = 2m$ for some $m \geq 1$.  
Since $\alpha_{\sigma} \in G(S^{2m}_{\Q},X)$, the map 
$\alpha_{\sigma} : S^{2m}_{\Q} \to X$ can be extended to the James reduced product space by Proposition \ref{pr:jamesredprodspext}, say, 
$$
\overline{\alpha}_{\sigma} : (S^{2m}_{\Q})_\infty \longrightarrow X, \ \overline{\alpha}_{\sigma} \in G((S^{2m}_{\Q})_\infty,X),
$$
where we know $(S^{2m}_{\Q})_\infty$ $\simeq$ $(S^{2m}_\infty)_{\Q}$ $\simeq$ $(\Omega \Sigma S^{2m})_{\Q}$ $\simeq$ $(\Omega S^{2m{+}1})_{\Q}$ $\simeq$ $\Omega (S^{2m+1}_{\Q})$ $\simeq$ $\Omega K(\Q,2m{+}1)$ $\simeq$ $K(\Q,2m)$.
Thus we have $\overline{\alpha}_{\sigma} \in G(K(\Q,2m),X)$.
Hence by Proposition \ref{pr:infpairingext}, there is a map $\overline{\alpha} : K(R,2m) = \bigoplus_{\lambda}K(\Q,2m) \to X$ extending $\alpha : S^{n}(R) \to X$.
Then we obtain $\beta\comp\overline{\alpha} \sim {\rm id}_{K(R,n)}$, since the identity map ${\rm id} : K(R,n) \to K(R,n)$ is the unique extension of $\iota^{n}_{R} : S^{n}(R) \to K(R,n)$, up to homotopy.
\par
Thus in either case, we obtain a map $\overline{\alpha} \in G(K(R,n),X)$ such that 
$$
\beta \comp \overline{\alpha} \sim {\rm id} : 
K(R,n) \longrightarrow K(R,n).
$$ 
Let $Y$ be the homotopy fibre of $\beta : X \to K(R,n)$.
Then by Theorem 1.3 of Dula and Gottlieb \cite{DuGo90}, we obtain 
$$
X \simeq Y \times \bigoplus_{\lambda\in\Lambda} K(\Q,n)
\simeq Y \times K(R,n).
$$
This completes the proof of the theorem.
\end{proof}  

\section{Applications}

A $0$-connected space $X$ is called a $T$-{\it space} if the fibration $\Omega X \to X^{S^1} \to X$ is trivial in the sense of fibre homotopy type (Aguad\'{e} \cite{Aguade87}). 
If $X$ is a $0$-connected Hopf space, then $X$ is a $T$-space. 
Aguad\'{e} showed that $1$-connected space $X$ of finite type is a rational $T$-space if and only if $X$ has the same rational homotopy type as a generalized Eilenberg-Mac{\,}Lane space, i.e., a product of (infinitely many) Eilenberg-Mac{\,}Lane spaces (Theorem 3.3 of \cite{Aguade87}).
Woo and Yoon  showed that a space $X$ is a $T$-space if and only if $G(\Sigma A, X) =[\Sigma A, X]$ for any space $A$ by Theorem 2.2 of \cite{WoYo95}.  So, it might be more appropriate to call such space a generalized Gottlieb space.
Then we have the following result by Theorem \ref{th:rathopfdec}.

\begin{thm}\label{th:ratTspsplit}
Let $R=\bigoplus_{\lambda \in \Lambda}\Q\,$ be  a finite or an infinite dimensional $\Q$-vector space.  Let $X$ be a $0$-connected $T$-space and  $R \subset \pi_n (X)$, $n\geq 2$.   
If $\overline{\rho} | R : R \to H_n(X)$ is an injection, where $R$ is represented as $R \subset [S_{\Q}^{n},X]=G(S_{\Q}^{n},X)$.  
Then  $X$ decomposes as 
$$
X\simeq Y \times \bigoplus_{\lambda\in\Lambda}K(\Q,n) \simeq Y \times K(R,n),
\quad \text{for a $T$-space $Y$.}
$$ 
\end{thm} 
\begin{proof}  By Theorem 2.11 of \cite{WoYo95} 
and Theorem \ref{th:rathopfdec},  
we have the result.  
\end{proof}  

Theorem \ref{th:ratTspsplit} implies the following result as a 
direct consequence.

\begin{cor}\label{cor:1connTspsplit} 
Let $n\geq 2$.  Let $R=\bigoplus_{\lambda \in \Lambda}\Q\,$ 
be a finite or an infinite dimensional $\Q$-vector space and assume that $R \subset \pi_n (X)$.  
If $X$ is $(n-1)$-connected $T$-space, then $X$ splits as 
$$
X\simeq Y \times \bigoplus_{\lambda\in\Lambda}K(\Q,n) \simeq Y \times K(R,n),
\quad \text{for a $T$-space $Y$.}
$$ 
\end{cor} 

A space $X$ is called a $G$-{\it space} if $G_n (X) = \pi_n(X)$ for all $n$ (cf.\ \cite{Gottlieb69a}).  As a special case of Theorem \ref{th:rathopfdec}, we have the following result for rational Gottlieb space. We remark that $ \pi_n(X_\Q)=G_n (X_\Q)$  implies $[S^{n}_{\Q}, X_{\Q}] = G (S^{n}_{\Q}, X_{\Q})$ for any $n$.

\begin{thm}\label{th:nm1connGspsplit}
Let $n\geq 2$.  Assume that a rational space $X_\Q$ is an $(n{-}1)$-connected Gottlieb space.   
If $\pi_n(X_\Q) \cong \bigoplus_{\lambda \in \Lambda}\Q\,$,  a finite or an infinite dimensional $\Q$-vector space, then $X_\Q$ decomposes as 
$$
X_\Q \simeq Y_\Q \times \bigoplus_{\lambda\in\Lambda}K(\Q,n) \simeq Y_\Q \times K(\pi_n(X_\Q),n),
$$ 
where $Y_\Q$ is an $n$-connected Gottlieb space.  
\end{thm} 

Theorem \ref{th:nm1connGspsplit} implies the following theorem (cf. \cite{Scheerer85}).  
For finite complexes or finite Postnikov pieces, it is known by Haslam \cite{Haslam71} and Mataga \cite{Mataga74}. 

\begin{thm}\label{th:HspTspGsp} 
If $X$ is a $1$-connected space, then the following  are equivalent$:$ \par
$(1)$  \ $X_\Q$ is a Gottlieb space. \par
$(2)$  \ $X_\Q$ is a $T$-space. \par
$(3)$  \ $X_\Q$ is a Hopf space. \par
$(4)$  \ $X_\Q$ has the homotopy type of a generalized Eilenberg-Mac{\,}Lane space.
\end{thm} 


\begin{cor}\label{cor:kinvariant}
Any $k$-invariant of a $1$-connected Gottlieb space is rationally trivial.
\end{cor}
We remark that Corollary \ref{cor:kinvariant} doesn't imply that a $k$-invariant of a $1$-connected Gottlieb space is of finite order.
Now, $H_{\ast}(K(\oplus_{\lambda}{\mathbb Q},2m{+}1);{\mathbb Q})$ is isomorphic to an exterior alebra and $H_{\ast}(K(\oplus_{\lambda}{\mathbb Q},2m);{\mathbb Q})$ is isomorphic to a polynomial algebra as Hopf algebras.
Thus we obtain a generalization of Theorem 3.2 of Borel \cite{Borel67}:


\begin{cor}\label{cor:genBorel}  
Let $X$ be a $1$-connected rational Gottlieb space.
Then $X_\Q$ is a Hopf space and the Hopf algebra $H^{\ast}(X;{\mathbb Q})$ is isomorphic (as an algebra) to the tensor product of the dual algebra of a polynomial algebra on even degree generators and the dual algebra of an exterior algebra on odd degree generators.
\end{cor}
We remark that $\pi_{q}(X){\otimes}{\mathbb Q}$ may be infinite dimensional for each $q \geq 1$, and hence $H_{q}(X;{\mathbb Q})$ and its dual $H^{q}(X;{\mathbb Q}) \cong \Hom(H_{q}(X;{\mathbb Q});\Q)$ may be distinct as $\Q$-modules for each $q \geq 1$.
For example, the dual of an exterior algebra on $\{\alpha_{\lambda}\}$ is not an exterior algebra on $\{\bar\alpha_{\lambda}\}$, in general, where $\bar\alpha_{\lambda}$ is the dual to $\alpha_{\lambda}$ (cf. \cite{MacLane:homology}).

\end{document}